\documentclass[12pt]{article}
\usepackage{amsmath}
\usepackage{amsfonts}
\usepackage{amssymb}

\def\thebibliograph#1#2{\section*{{\normalsize \bf #2}}\list
   {[\arabic{enumi}]}{\settowidth\labelwidth{[#1]}\leftmargin\labelwidth
     \advance\leftmargin\labelsep
     \usecounter{enumi}}
     \def\newblock{\hskip .11em plus .33em minus -.07em}
     \sloppy
     \sfcode`\.=1000\relax}

\newtheorem{theorem}{Theorem}
\newtheorem{prop}{Proposition}
\newtheorem{definition}{Definition}

\newtheorem{lemma}{Lemma}
\newtheorem{remark}{Remark}
\input{tcilatex}
\begin{document}

\title{On the duality of variable Triebel-Lizorkin spaces}
\author{ Douadi Drihem \ \thanks{%
M'sila University, Department of Mathematics, Laboratory of Functional
Analysis and Geometry of Spaces , P.O. Box 166, M'sila 28000, Algeria,
e-mail: \texttt{\ douadidr@yahoo.fr}}}
\date{\today }
\maketitle

\begin{abstract}
The aim of this paper is to prove duality of Triebel-Lizorkin spaces $%
F_{1,q\left( \cdot \right) }^{\alpha \left( \cdot \right) }$. First, we
prove the duality of associated sequence spaces. Then from the so-called $%
\varphi $-transform characterization in the sense of Frazier and Jawerth, we
deduce the main result of this paper.\vskip5pt

\textit{MSC 2010\/}: 46B10, 46E35.

\textit{Key Words and Phrases}: Besov-type space, Triebel-Lizorkin spaces,
Duality, Variable exponent.
\end{abstract}

\section{Introduction}

In recent years, there has been growing interest in generalizing classical
spaces such as Lebesgue, Sobolev spaces, Besov spaces, Triebel-Lizorkin
spaces to the case with either variable integrability or variable
smoothness. The motivation for the increasing interest in such spaces comes
not only from theoretical purposes, but also from applications to fluid
dynamics \cite{Ru00}, image restoration \cite{CLR06} and PDE with
non-standard growth conditions. Leopold [17, 18, 19, 20] and Leopold \&
Schrohe [21] studied pseudo-differential operators, they introduced related
Besov spaces with variable smoothness $B_{p,p}^{\alpha (\cdot )}$. Function
spaces of variable smoothness have recently been studied by Besov [2, 3, 4].
Along a different line of study, J.-S. Xu \cite{Xu08}, \cite{Xu09} has
studied Besov spaces with variable $p$, but fixed $q$ and $\alpha $.\vskip5pt

Besov spaces of variable smoothness and integrability, $B_{p(\cdot ),q(\cdot
)}^{\alpha (\cdot )}$, initially appeared in the paper of A. Almeida and P. H%
\"{a}st\"{o} \cite{AH}. Several basic properties were established, such as
the Fourier analytical characterisation and Sobolev embeddings. When $%
p,q,\alpha $ are constants they coincide with the usual function spaces $%
B_{p,q}^{s}$. \vskip5pt

Variable Besov-type spaces have been introduced in \cite{D5} and \cite{D6},
where their basic properties are given, such as the Sobolev type embeddings
and \ that under some conditions these spaces are just the variable Besov
spaces. For constant exponents, these spaces unify and generalize many
classical function spaces including Besov spaces, Besov-Morrey spaces (see,
for example, \cite[Corollary 3.3]{WYY}). Independently, D. Yang, C. Zhuo and
W. Yuan, \cite{YZW15} studied these function spaces where several properties
are obtained such as atomic decomposition and the boundedness of trace
operator, see \cite{YZW151} for further properties of variable
Triebel-Lizorkin-type spaces. Also, A. I. Tyulenev \cite{Ty151}, \cite{Ty152}
has studied some new function spaces of variable smoothness.
Triebel-Lizorkin spaces with variable exponents $F_{p(\cdot ),q(\cdot
)}^{\alpha (\cdot )}$ were introduced by \cite{DHR}. They proved a
discretization by the so called $\varphi $-transform. Also atomic and
molecular decomposition of these function spaces are obtained and used it to
derive trace results. The Sobolev embedding of these function spaces was
proved by J. Vyb\'{\i}ral, \cite{V}. Some properties of these function
spaces such as local means characterizations and characterizations by ball
means of differences can be found in \cite{KV121} and \cite{KV122}. When $%
\alpha ,p,q$ are constants they coincide with the usual function spaces $%
F_{p,q}^{\alpha }$.

It is well-known that duality is an important concept when we study function
spaces. It applied to real interpolation and embeddings. See \cite{T73}, for
the duality of the usual Besov spaces $B_{p,q}^{\alpha }$ and
Triebel-Lizorkin spaces $F_{p,q}^{\alpha }$. M. Izuki and T. Noi \cite{IN14}
have obtained the duality of $B_{p(\cdot ),q(\cdot )}^{\alpha (\cdot )}$ and 
$F_{p(\cdot ),q(\cdot )}^{\alpha (\cdot )}$, for $1<p^{-}\leq p^{+}<\infty $
and $1<q^{-}\leq q^{+}<\infty $, see \cite{N14} for the duality of $%
B_{p(\cdot ),q}^{\alpha }$ and $F_{p(\cdot ),q}^{\alpha }$ spaces with
variable $p$, but fixed $q$ and $\alpha $. In the present paper we obtain
the duality of variable Triebel-Lizorkin spaces $F_{1,q(\cdot )}^{\alpha
(\cdot )}$.

\section{Preliminaries}

As usual, we denote by $\mathbb{R}^{n}$ the $n$-dimensional real Euclidean
space, $\mathbb{N}$ the collection of all natural numbers and $\mathbb{N}%
_{0}=\mathbb{N}\cup \{0\}$. The letter $\mathbb{Z}$ stands for the set of
all integer numbers. The expression $f\lesssim g$ means that $f\leq c\,g$
for some independent constant $c$ (and non-negative functions $f$ and $g$),
and $f\approx g$ means $f\lesssim g\lesssim f$.\vskip5pt

By supp $f$ we denote the support of the function $f$ , i.e., the closure of
its non-zero set. If $E\subset {\mathbb{R}^{n}}$ is a measurable set, then $%
|E|$ stands for the (Lebesgue) measure of $E$ and $\chi _{E}$ denotes its
characteristic function.\vskip5pt

The symbol $\mathcal{S}(\mathbb{R}^{n})$ is used in place of the set of all
Schwartz functions on $\mathbb{R}^{n}$. We denote by $\mathcal{S}^{\prime }(%
\mathbb{R}^{n})$ the dual space of all tempered distributions on $\mathbb{R}%
^{n}$. We define the Fourier transform of a function $f\in \mathcal{S}(%
\mathbb{R}^{n})$ by $\mathcal{F}(f)(\xi )=(2\pi )^{-n/2}\int_{\mathbb{R}%
^{n}}e^{-ix\cdot \xi }f(x)dx$. Its inverse is denoted by $\mathcal{F}^{-1}f$%
. Both $\mathcal{F}$ and $\mathcal{F}^{-1}$ are extended to the dual
Schwartz space $\mathcal{S}^{\prime }(\mathbb{R}^{n})$ in the usual way.%
\vskip5pt

For $v\in \mathbb{Z}$ and $m=(m_{1},...,m_{n})\in \mathbb{Z}^{n}$, let $%
Q_{v,m}$ be the dyadic cube in $\mathbb{R}^{n}$, $Q_{v,m}=%
\{(x_{1},...,x_{n}):m_{i}\leq 2^{v}x_{i}<m_{i}+1,i=1,2,...,n\}$. For the
collection of all such cubes we use $\mathcal{Q}:=\{Q_{v,m}:v\in \mathbb{Z}%
,m\in \mathbb{Z}^{n}\}$. For each cube $Q$, we denote by $x_{Q_{v,m}}$ the
lower left-corner $2^{-v}m$ of $Q=Q_{v,m}$, its side length by $l(Q)$ and
for $r>0$, we denote by $rQ$ the cube concentric with $Q$ having the side
length $rl(Q)$. Furthermore, we put $v_{Q}=-\log _{2}l(Q)$, $v_{Q}^{+}=\max
(v_{Q},0)$ and $\chi _{v,m}=\chi _{Q_{v,m}}$.\vskip5pt

For $v\in \mathbb{Z}$, $\varphi \in \mathcal{S}(\mathbb{R}^{n})$ and $x\in 
\mathbb{R}^{n}$, we set $\widetilde{\varphi }(x):=\overline{\varphi (-x)}$, $%
\varphi _{v}(x):=2^{vn}\varphi (2^{v}x)$, and%
\begin{equation*}
\varphi _{v,m}(x):=2^{vn/2}\varphi (2^{v}x-m)=|Q_{v,m}|^{1/2}\varphi
_{v}(x-x_{Q_{v,m}})\quad \text{if\quad }Q=Q_{v,m}.
\end{equation*}

By $c$ we denote generic positive constants, which may have different values
at different occurrences. Although the exact values of the constants are
usually irrelevant for our purposes, sometimes we emphasize their dependence
on certain parameters (e.g. $c(p)$ means that $c$ depends on $p$, etc.).
Further notation will be properly introduced whenever needed.

The variable exponents that we consider are always measurable functions $p$
on $\mathbb{R}^{n}$ with range in $[c,\infty \lbrack $ for some $c>0$. We
denote the set of such functions by $\mathcal{P}_{0}$. The subset of
variable exponents with range $[1,\infty \lbrack $ is denoted by $\mathcal{P}
$. We use the standard notation $p^{-}:=\underset{x\in \mathbb{R}^{n}}{\text{%
ess-inf}}$ $p(x)$,$\quad p^{+}:=\underset{x\in \mathbb{R}^{n}}{\text{ess-sup 
}}p(x)$.

The variable exponent modular is defined by $\varrho _{p(\cdot )}(f):=\int_{%
\mathbb{R}^{n}}\varrho _{p(x)}(\left\vert f(x)\right\vert )dx$, where $%
\varrho _{p}(t)=t^{p}$. The variable exponent Lebesgue space $L^{p(\cdot )}$%
\ consists of measurable functions $f$ on $\mathbb{R}^{n}$ such that $%
\varrho _{p(\cdot )}(\lambda f)<\infty $ for some $\lambda >0$. We define
the Luxemburg (quasi)-norm on this space by the formula $\left\Vert
f\right\Vert _{p(\cdot )}:=\inf \big\{\lambda >0:\varrho _{p(\cdot )}\Big(%
\frac{f}{\lambda }\Big)\leq 1\big\}$. A useful property is that $\left\Vert
f\right\Vert _{p(\cdot )}\leq 1$ if and only if $\varrho _{p(\cdot )}(f)\leq
1$, see \cite{DHHR}, Lemma 3.2.4. For variable exponents, H\"{o}lder's
inequality takes the form $\Vert fg\Vert _{s(\cdot )}\lesssim \,\Vert f\Vert
_{{p(\cdot )}}\Vert g\Vert _{q(\cdot )}$ where $s$ is defined pointwise by $%
\tfrac{1}{s(x)}=\tfrac{1}{p(x)}+\tfrac{1}{q(x)}$. Often we use the
particular case $s(x):=1$ corresponding to the situation when $q=p^{\prime }$
is the conjugate exponent of $p$.

Let $p,q\in \mathcal{P}_{0}$. The mixed Lebesgue-sequence space $\ell
^{q(\cdot )}(L^{p(\cdot )})$ is defined on sequences of $L^{p(\cdot )}$%
-functions by the modular%
\begin{equation*}
\varrho _{\ell ^{q(\cdot )}(L^{p\left( \cdot \right)
})}((f_{v})_{v}):=\sum\limits_{v}\inf \Big\{\lambda _{v}>0:\varrho _{p(\cdot
)}\Big(\frac{f_{v}}{\lambda _{v}^{1/q(\cdot )}}\Big)\leq 1\Big\}.
\end{equation*}%
The (quasi)-norm is defined from this as usual:%
\begin{equation}
\left\Vert \left( f_{v}\right) _{v}\right\Vert _{\ell ^{q(\cdot
)}(L^{p\left( \cdot \right) })}:=\inf \Big\{\mu >0:\varrho _{\ell ^{q(\cdot
)}(L^{p(\cdot )})}\Big(\frac{1}{\mu }(f_{v})_{v}\Big)\leq 1\Big\}.
\label{mixed-norm}
\end{equation}%
If $q^{+}<\infty $, then we can replace $\mathrm{\eqref{mixed-norm}}$ by the
simpler expression $\varrho _{\ell ^{q(\cdot )}(L^{p(\cdot
)})}((f_{v})_{v}):=\sum\limits_{v}\left\Vert |f_{v}|^{q(\cdot )}\right\Vert
_{\frac{p(\cdot )}{q(\cdot )}}$. Furthermore, if $p$ and $q$ are constants,
then $\ell ^{q(\cdot )}(L^{p(\cdot )})=\ell ^{q}(L^{p})$. The case $%
p:=\infty $ can be included by replacing the last modular by $\varrho _{\ell
^{q(\cdot )}(L^{\infty })}((f_{v})_{v}):=\sum\limits_{v}\big\|\left\vert
f_{v}\right\vert ^{q(\cdot )}\big\|_{\infty }$. Let $p$, $q\in \mathcal{P}%
_{0}$. Then $\varrho _{\ell ^{q(\cdot )}(L^{p(\cdot )})}$ is continuous if $%
p^{+}<\infty $ and $q^{+}<\infty $, see\ \cite{AH}.

It is known, cf. \cite{AH} and \cite{KV121}, that $\ell ^{q(\cdot
)}(L^{p(\cdot )})$ is a norm if $q(\cdot )\geq 1$ is constant almost
everywhere (a.e.) on $\mathbb{R}^{n}$ and $p(\cdot )\geq 1$, or if $\frac{1}{%
p(x)}+\frac{1}{q(x)}\leq 1$ a.e. on $\mathbb{R}^{n}$, or if $1\leq q(x)\leq
p(x)<\infty $ a.e. on $\mathbb{R}^{n}$.

We state also the definition of the space $L^{p(\cdot )}(\ell ^{q(\cdot )})$
which is much more intuitive then the definition of $\ell ^{q(\cdot
)}(L^{p(\cdot )})$. One just takes the $\ell ^{q(x)}$ norm of $%
(f_{v}(x))_{v} $ for every $x\in \mathbb{R}^{n}$ and then the $L^{p(\cdot )}$%
-norm with respect to $x\in \mathbb{R}^{n}$, i.e.%
\begin{equation*}
\big\|\left( f_{v}\right) _{v\geq 0}\big\|_{L^{p(\cdot )}(\ell ^{q(\cdot
)})}:=\big\|\big\|\left( f_{v}(x)\right) _{v\geq 0}\big\|_{\ell ^{q(x)}}%
\big\|_{p(\cdot )}.
\end{equation*}%
It is easy to show that $L^{p(\cdot )}(\ell ^{q(\cdot )})$\ is always a
quasi-normed space\ and it is a normed space, if $\min (p(x),q(x))\geq 1$\
holds point-wise.

We say that $g:\mathbb{R}^{n}\rightarrow \mathbb{R}$ is \textit{locally }log%
\textit{-H\"{o}lder continuous}, abbreviated $g\in C_{\text{loc}}^{\log }$,
if there exists $c_{\log }(g)>0$ such that%
\begin{equation}
\left\vert g(x)-g(y)\right\vert \leq \frac{c_{\log }(g)}{\log
(e+1/\left\vert x-y\right\vert )}  \label{lo-log-Holder}
\end{equation}%
for all $x,y\in \mathbb{R}^{n}$. We say that $g$ satisfies the log\textit{-H%
\"{o}lder decay condition}, if there exists $g_{\infty }\in \mathbb{R}$ and
a constant $c_{\log }>0$ such that%
\begin{equation*}
\left\vert g(x)-g_{\infty }\right\vert \leq \frac{c_{\log }}{\log
(e+\left\vert x\right\vert )}
\end{equation*}%
for all $x\in \mathbb{R}^{n}$. We say that $g$ is \textit{globally}-log%
\textit{-H\"{o}lder continuous}, abbreviated $g\in C^{\log }$, if it is%
\textit{\ }locally log-H\"{o}lder continuous and satisfies the log-H\"{o}%
lder decay\textit{\ }condition.\textit{\ }The constants $c_{\log }(g)$ and $%
c_{\log }$ are called the \textit{locally }log\textit{-H\"{o}lder constant }%
and the log\textit{-H\"{o}lder decay constant}, respectively\textit{.} We
note that all functions $g\in C_{\text{loc}}^{\log }$ always belong to $%
L^{\infty }$.\vskip5pt

We define the following class of variable exponents $\mathcal{P}^{\mathrm{log%
}}:=\big\{p\in \mathcal{P}:\frac{1}{p}$ is globally-log-H\"{o}lder continuous%
$\big\}$, were introduced in $\mathrm{\cite[Section \ 2]{DHHMS}}$. We define 
$1/p_{\infty }:=\lim_{|x|\rightarrow \infty }1/p(x)$\ and we use the
convention $\frac{1}{\infty }=0$. Note that although $\frac{1}{p}$ is
bounded, the variable exponent $p$ itself can be unbounded. If $p\in 
\mathcal{P}^{\mathrm{log}}$, then the convolution with a radially decreasing 
$L^{1}$-function is bounded on $L^{p(\cdot )}$: $\Vert \varphi \ast f\Vert _{%
{p(\cdot )}}\leq c\Vert \varphi \Vert _{{1}}\Vert f\Vert _{{p(\cdot )}}$. %
\vskip5pt

It is known that for $p\in \mathcal{P}^{\mathrm{log}}$ we have%
\begin{equation}
\Vert \chi _{B}\Vert _{{p(\cdot )}}\Vert \chi _{B}\Vert _{{p}^{\prime }{%
(\cdot )}}\approx |B|.  \label{DHHR}
\end{equation}%
Also,%
\begin{equation}
\Vert \chi _{B}\Vert _{{p(\cdot )}}\approx |B|^{\frac{1}{p(x)}},\quad x\in B
\label{DHHR1}
\end{equation}%
for small balls $B\subset {\mathbb{R}^{n}}$ ($|B|\leq 2^{n}$), with
constants only depending on the $\log $-H\"{o}lder constant of $p$ (see, for
example, \cite[Section 4.5]{DHHR}). These properties are hold if $p\in 
\mathcal{P}_{0}^{\mathrm{log}}$, since $\Vert \chi _{B}\Vert _{{p(\cdot )}%
}=\Vert \chi _{B}\Vert _{{p(\cdot )/a}}^{1/a}$ and $\frac{p}{a}{\in }%
\mathcal{P}^{\mathrm{log}}$ if $p^{-}\geq a.$\vskip5pt

Recall that $\eta _{v,N}(x):=2^{nv}(1+2^{v}\left\vert x\right\vert )^{-N}$,
for any $x\in \mathbb{R}^{n}$, $v\in \mathbb{N}_{0}$ and $N>0$. Note that $%
\eta _{v,N}\in L^{1}$ when $N>n$ and that $\left\Vert \eta _{v,N}\right\Vert
_{1}=c_{N}$ is independent of $v$. We introduce the abbreviations%
\begin{equation*}
\left\Vert \left( f_{v}\right) _{v}\right\Vert _{\ell ^{q(\cdot
)}(L_{p(\cdot )}^{p(\cdot )})}:=\sup_{\{P\in \mathcal{Q},|P|\leq 1\}}\Big\|%
\Big(\frac{f_{v}}{|P|^{1/p(\cdot )}}\chi _{P}\Big)_{v\geq v_{P}}\Big\|_{\ell
^{q(\cdot )}(L^{p(\cdot )})}.
\end{equation*}

The following lemma is the $\ell ^{q(\cdot )}(L_{p(\cdot )}^{p(\cdot )})$%
-version of Lemma 4.7 from A. Almeida and P. H\"{a}st\"{o} \cite{AH} (we use
it, since the maximal operator is in general not bounded on $\ell ^{q(\cdot
)}(L^{p(\cdot )})$, see \cite[Example 4.1]{AH}).

\begin{lemma}
\label{Alm-Hastolemma1}Let $p\in \mathcal{P}^{\log }$, $q\in \mathcal{P}%
_{0}^{\log }$ with $0<q^{-}\leq q^{+}<\infty $ and $p^{-}>1$. For $%
N>2n+c_{\log }(1/p)+c_{\log }(1/q)$, there exists $c>0$ such that%
\begin{equation*}
\left\Vert (\eta _{v,N}\ast f_{v})_{v}\right\Vert _{\ell ^{q(\cdot
)}(L_{p(\cdot )}^{p(\cdot )})}\leq c\left\Vert (f_{v})_{v}\right\Vert _{\ell
^{q(\cdot )}(L_{p(\cdot )}^{p(\cdot )})}.
\end{equation*}
\end{lemma}

The arguments in \cite[Lemma 2.12]{D5}, are true to prove this property, in
view of the fact that $\left\Vert \chi _{P}\right\Vert _{p(\cdot )}\approx
|P|^{1/p(\cdot )}$, since $|P|\leq 1$. The proof of the following lemma is
given in \cite[Theorem 3.2]{DHR}$\mathrm{.}$

\begin{lemma}
\label{Alm-Hastolemma2}Let $p,q\in \mathcal{P}^{\log }$ with $1<p^{-}\leq
p^{+}<\infty $ and $1<q^{-}\leq q^{+}<\infty $. For $N>n$, there exists $c>0$
such that%
\begin{equation*}
\left\Vert \left( \eta _{v,N}\ast f_{v}\right) _{v}\right\Vert _{L^{p(\cdot
)}(\ell ^{q(\cdot )})}\leq c\left\Vert \left( f_{v}\right) _{v}\right\Vert
_{L^{p(\cdot )}(\ell ^{q(\cdot )})}.
\end{equation*}
\end{lemma}

\section{Spaces of variable smoothness and integrability}

In this section we recall the definition of the spaces $\widetilde{B}%
_{p(\cdot ),q(\cdot )}^{\alpha (\cdot ),p(\cdot )}$ and $F_{p(\cdot
),q(\cdot )}^{\alpha (\cdot )}$\ as given in \cite{D6} and \cite{DHR}. Let $%
\Psi $\ be a function\ in $\mathcal{S}(\mathbb{R}^{n})$\ satisfying $0\leq
\Psi (x)\leq 1$ for all $x$, $\Psi (x)=1$\ for\ $\left\vert x\right\vert
\leq 1$\ and\ $\Psi (x)=0$\ for\ $\left\vert x\right\vert \geq 2$.\ We put $%
\mathcal{F}\varphi _{0}(x)=\Psi (x)$, $\mathcal{F}\varphi (x)=\Psi (\frac{x}{%
2})-\Psi (x)$\ and $\mathcal{F}\varphi _{v}(x)=\mathcal{F}\varphi
(2^{-v+1}x) $ for $v=1,2,3,....$ Then $\{\mathcal{F}\varphi _{v}\}_{v\in 
\mathbb{N}_{0}}$\ is a resolution of unity, $\sum_{v=0}^{\infty }\mathcal{F}%
\varphi _{v}(x)=1 $ for all $x\in \mathbb{R}^{n}$.\ Thus we obtain the
Littlewood-Paley decomposition 
\begin{equation}
f=\sum_{v=0}^{\infty }\varphi _{v}\ast f  \label{L-P-decomposition}
\end{equation}%
of all $f\in \mathcal{S}^{\prime }(\mathbb{R}^{n})$ $($convergence in $%
\mathcal{S}^{\prime }(\mathbb{R}^{n}))$.

Now, we define the spaces under consideration.

\begin{definition}
\label{B-F-def}Let $\{\mathcal{F}\varphi _{v}\}_{v\in \mathbb{N}_{0}}$ be a
resolution of unity, $\alpha :\mathbb{R}^{n}\rightarrow \mathbb{R}$ and $%
p,q\in \mathcal{P}_{0}$.

$\mathrm{(i)}$ The Besov-type space $\widetilde{B}_{p(\cdot ),q(\cdot
)}^{\alpha (\cdot ),p(\cdot )}$\ is the collection of all $f\in \mathcal{S}%
^{\prime }(\mathbb{R}^{n})$\ such that 
\begin{equation}
\left\Vert f\right\Vert _{\widetilde{B}_{p(\cdot ),q(\cdot )}^{\alpha (\cdot
),p(\cdot )}}:=\sup_{P\in \mathcal{Q}}\Big\|\Big(\frac{2^{v\alpha \left(
\cdot \right) }\varphi _{v}\ast f}{|P|^{1/p(\cdot )}}\chi _{P}\Big)_{v\geq
v_{P}^{+}}\Big\|_{\ell ^{q(\cdot )}(L^{p(\cdot )})}<\infty ,  \label{B-def}
\end{equation}%
$\mathrm{(ii)}$ Let $0<p^{-}\leq p^{+}<\infty $ and $0<q^{-}\leq
q^{+}<\infty $. The Triebel-Lizorkin space $F_{p(\cdot ),q(\cdot )}^{\alpha
(\cdot )}$\ is the collection of all $f\in \mathcal{S}^{\prime }(\mathbb{R}%
^{n})$\ such that 
\begin{equation}
\left\Vert f\right\Vert _{F_{p(\cdot ),q(\cdot )}^{\alpha (\cdot
)}}:=\left\Vert \left( 2^{\alpha \left( \cdot \right) }\varphi _{v}\ast
f\right) _{v\geq 0}\right\Vert _{L^{p(\cdot )}(\ell ^{q(\cdot )})}<\infty .
\label{F-def}
\end{equation}
\end{definition}

The definition of the spaces $\widetilde{B}_{p(\cdot ),q(\cdot )}^{\alpha
(\cdot ),p(\cdot )}$\ and $F_{p(\cdot ),q(\cdot )}^{\alpha (\cdot )}$ is
independent of the chosen resolution of unity $\mathrm{%
\eqref{L-P-decomposition}}$\textrm{\ }if $\alpha \in C_{\mathrm{loc}}^{\log
} $, $p,q\in \mathcal{P}_{0}^{\log }$ and $0<q^{+}<\infty $, ($0<p^{-}\leq
p^{+}<\infty $ in the $F$ case) and that different choices yield equivalent
quasi-norms. Using the system $\{\varphi _{v}\}_{v\in \mathbb{N}_{0}}$ we
can define the norm%
\begin{equation*}
\left\Vert f\right\Vert _{B_{p,q}^{\alpha ,\tau }}:=\sup_{P\in \mathcal{Q}}%
\frac{1}{\left\vert P\right\vert ^{\tau }}\Big(\sum\limits_{v=v_{P}^{+}}^{%
\infty }2^{v\alpha q}\left\Vert (\varphi _{v}\ast f)\chi _{P}\right\Vert
_{p}^{q}\Big)^{1/q}
\end{equation*}%
for constants $\alpha $ and $p,q\in (0,\infty ]$. The Besov-type space $%
B_{p,q}^{\alpha ,\tau }$ consist of all distributions $f\in \mathcal{S}%
^{\prime }(\mathbb{R}^{n})$ for which $\left\Vert f\right\Vert
_{B_{p,q}^{\alpha ,\tau }}<\infty $. It is well-known that these spaces do
not depend on the choice of the system $\{\mathcal{F}\varphi _{v}\}_{v\in 
\mathbb{N}_{0}}$ (up to equivalence of quasinorms). Further details on the
classical theory of these spaces can be found in \cite{D2} and \cite{WYY};
see also \cite{D4} for recent developments. One recognizes immediately that
if $\alpha $, $p$ and $q$ are constants, then $\widetilde{B}_{p(\cdot
),q(\cdot )}^{\alpha (\cdot ),p(\cdot )}=B_{p,q}^{\alpha ,1/p}$ and $%
F_{p(\cdot ),q(\cdot )}^{\alpha (\cdot )}=F_{p,q}^{\alpha }$ is the
classical Triebel-Lizorkin spaces, see \cite{T82}, \cite{T92} and \cite{RS96}
for the history of these spaces. When, $q:=\infty $\ the Besov-type space $%
\widetilde{B}_{p(\cdot ),\infty }^{\alpha (\cdot ),p(\cdot )}$\ consist of
all distributions $f\in \mathcal{S}^{\prime }(\mathbb{R}^{n})$\ such that 
\begin{equation*}
\sup_{P\in \mathcal{Q},v\geq v_{P}^{+}}\Big\|\frac{2^{v\alpha \left( \cdot
\right) }\varphi _{v}\ast f}{|P|^{1/p(\cdot )}}\chi _{P}\Big\|_{p(\cdot
)}<\infty .
\end{equation*}%
If we replace dyadic cubes $P$ in Definition \ref{B-F-def} by arbitrary
cubes $P$, we then obtain equivalent quasi-norms. It is clear that if $%
\alpha $ and $p$ are constants, then $\widetilde{B}_{p(\cdot ),p(\cdot
)}^{\alpha (\cdot ),p(\cdot )}=F_{\infty ,p}^{\alpha }$, see \cite{FJ90} for
the properties of $F_{\infty ,p}^{\alpha }$. We refer to the papers $\mathrm{%
\cite{D5}}$\textrm{, }$\mathrm{\cite{D6}}$ and $\mathrm{\cite{YZW15}}$ 
\textrm{,} where various results on variable Besov-type spaces were obtained.

\begin{lemma}
\label{new-equinorm copy(1)}A tempered distribution $f$ belongs to $%
\widetilde{B}_{p(\cdot ),q(\cdot )}^{\alpha (\cdot ),p(\cdot )}$ if and only
if,%
\begin{equation*}
\left\Vert f\right\Vert _{\widetilde{B}_{p(\cdot ),q(\cdot )}^{\alpha (\cdot
),p(\cdot )}}^{\#}:=\sup_{\{P\in \mathcal{Q},|P|\leq 1\}}\Big\|\Big(\frac{%
2^{v\alpha \left( \cdot \right) }\varphi _{v}\ast f}{|P|^{1/p(\cdot )}}\chi
_{P}\Big)_{v\geq v_{P}}\Big\|_{\ell ^{q(\cdot )}(L^{p(\cdot )})}<\infty .
\end{equation*}%
Furthermore, the quasi-norms $\left\Vert f\right\Vert _{\widetilde{B}%
_{p(\cdot ),q(\cdot )}^{\alpha (\cdot ),p(\cdot )}}$ and $\left\Vert
f\right\Vert _{\widetilde{B}_{p(\cdot ),q(\cdot )}^{\alpha (\cdot ),p(\cdot
)}}^{\#}$ are equivalent.
\end{lemma}

The proof this lemma is given in \cite{D6}. One of the key tools to prove
the duality of Triebel-Lizorkin spaces $F_{1,q\left( \cdot \right) }^{\alpha
\left( \cdot \right) }$ is to transfer the problem from function spaces to
their corresponding sequence spaces.

\begin{definition}
\label{sequence-space}Let $p,q\in \mathcal{P}_{0}$ and let $\alpha :\mathbb{R%
}^{n}\rightarrow \mathbb{R}$. Then for all complex valued sequences $\lambda
=\{\lambda _{v,m}\in \mathbb{C}\}_{v\in \mathbb{N}_{0},m\in \mathbb{Z}^{n}}$
we define%
\begin{equation*}
\widetilde{b}_{p\left( \cdot \right) ,q\left( \cdot \right) }^{\alpha \left(
\cdot \right) ,p(\cdot )}:=\Big\{\lambda :\left\Vert \lambda \right\Vert _{%
\widetilde{b}_{p\left( \cdot \right) ,q\left( \cdot \right) }^{\alpha \left(
\cdot \right) ,p(\cdot )}}<\infty \Big\}
\end{equation*}%
where%
\begin{equation*}
\left\Vert \lambda \right\Vert _{\widetilde{b}_{p\left( \cdot \right)
,q\left( \cdot \right) }^{\alpha \left( \cdot \right) ,p(\cdot
)}}:=\sup_{P\in \mathcal{Q}}\Big\|\Big(\frac{\sum\limits_{m\in \mathbb{Z}%
^{n}}2^{v(\alpha \left( \cdot \right) +n/2)}\lambda _{v,m}\chi _{v,m}}{%
|P|^{1/p(\cdot )}}\chi _{P}\Big)_{v\geq v_{P}^{+}}\Big\|_{\ell ^{q(\cdot
)}(L^{p(\cdot )})}
\end{equation*}
\end{definition}

If we replace dyadic cubes $P$ by arbitrary balls $B_{J}$ of $\mathbb{R}^{n}$
with $J\in \mathbb{Z}$, we then obtain equivalent quasi-norms, where the
supremum is taken over all $J\in \mathbb{Z}$\ and all balls $B_{J}$\ of $%
\mathbb{R}^{n}$.

\begin{remark}
Let $\alpha \in C_{\mathrm{loc}}^{\log }$, $p$, $q\in \mathcal{P}_{0}^{\log
} $ and $0<q^{+}<\infty $. Then 
\begin{equation*}
\left\Vert f\right\Vert _{\widetilde{B}_{q(\cdot ),q(\cdot )}^{\alpha (\cdot
),q(\cdot )}}\approx \sup_{\{P\in \mathcal{Q},|P|\leq 1\}}\Big\|\Big(\frac{%
2^{v\alpha \left( \cdot \right) }\varphi _{v}\ast f}{|P|^{1/q(\cdot )}}\chi
_{P}\Big)_{v\geq v_{P}}\Big\|_{L^{q(\cdot )}(\ell ^{q(\cdot )})}
\end{equation*}%
for any $f\in \widetilde{B}_{q(\cdot ),q(\cdot )}^{\alpha (\cdot ),q(\cdot
)} $. For all complex valued sequences $\lambda =\{\lambda _{v,m}\in \mathbb{%
C}:v\in \mathbb{N}_{0},m\in \mathbb{Z}^{n}\}\in \widetilde{b}_{q\left( \cdot
\right) ,q\left( \cdot \right) }^{\alpha \left( \cdot \right) ,q(\cdot )}$
we have%
\begin{equation}
\left\Vert \lambda \right\Vert _{\widetilde{b}_{q\left( \cdot \right)
,q\left( \cdot \right) }^{\alpha \left( \cdot \right) ,q(\cdot )}}\approx
\sup_{\{P\in \mathcal{Q},|P|\leq 1\}}\Big\|\Big(\frac{\sum\limits_{m\in 
\mathbb{Z}^{n}}2^{v(\alpha \left( \cdot \right) +n/2)}\lambda _{v,m}\chi
_{v,m}}{|P|^{1/q(\cdot )}}\chi _{P}\Big)_{v\geq v_{P}}\Big\|_{L^{q(\cdot
)}(\ell ^{q(\cdot )})}  \label{equiv-norm2}
\end{equation}%
for any $\lambda \in \widetilde{b}_{q\left( \cdot \right) ,q\left( \cdot
\right) }^{\alpha \left( \cdot \right) ,q(\cdot )}$, see $\mathrm{\cite{D6}}$%
.
\end{remark}

Let $\Phi $ and $\varphi $ satisfy%
\begin{equation}
\text{supp}\mathcal{F}\Phi \subset \overline{B(0,2)}\text{ and }|\mathcal{F}%
\Phi (\xi )|\geq c\text{ if }|\xi |\leq \frac{5}{3}  \label{Ass1}
\end{equation}%
and 
\begin{equation}
\text{supp}\mathcal{F}\varphi \subset \overline{B(0,2)}\backslash B(0,1/2)%
\text{ and }|\mathcal{F}\varphi (\xi )|\geq c\text{ if }\frac{3}{5}\leq |\xi
|\leq \frac{5}{3}  \label{Ass2}
\end{equation}%
where $c>0$. It easy to see that $\int_{\mathbb{R}^{n}}x^{\gamma }\varphi
(x)dx=0$ for all multi-indices $\gamma \in \mathbb{N}_{0}^{n}$. By \cite[pp.
130--131]{FJ90}, there exist \ functions $\Psi \in \mathcal{S}(\mathbb{R}%
^{n})$ satisfying $\mathrm{\eqref{Ass1}}$ and $\psi \in \mathcal{S}(\mathbb{R%
}^{n})$ satisfying $\mathrm{\eqref{Ass2}}$ such that for all $\xi \in 
\mathbb{R}^{n}$%
\begin{equation}
\mathcal{F}\widetilde{\Phi }(\xi )\mathcal{F}\Psi (\xi )+\sum_{j=1}^{\infty }%
\mathcal{F}\widetilde{\varphi }(2^{-j}\xi )\mathcal{F}\psi (2^{-j}\xi
)=1,\quad \xi \in \mathbb{R}^{n}.  \label{Ass4}
\end{equation}

Furthermore, we have the following identity for all $f\in \mathcal{S}%
^{\prime }(\mathbb{R}^{n})$; see \cite[(12.4)]{FJ90}%
\begin{eqnarray*}
f &=&\Psi \ast \widetilde{\Phi }\ast f+\sum_{v=1}^{\infty }\psi _{v}\ast 
\widetilde{\varphi }_{v}\ast f \\
&=&\sum_{m\in \mathbb{Z}^{n}}\widetilde{\Phi }\ast f(m)\Psi (\cdot
-m)+\sum_{v=1}^{\infty }2^{-vn}\sum_{m\in \mathbb{Z}^{n}}\widetilde{\varphi }%
_{v}\ast f(2^{-v}m)\psi _{v}(\cdot -2^{-v}m).
\end{eqnarray*}%
Recall that the $\varphi $-transform $S_{\varphi }$ is defined by setting $%
(S_{\varphi })_{0,m}=\langle f,\Phi _{m}\rangle $ where $\Phi _{m}(x)=\Phi
(x-m)$ and $(S_{\varphi })_{v,m}=\langle f,\varphi _{v,m}\rangle $ where $%
\varphi _{v,m}(x)=2^{vn/2}\varphi (2^{v}x-m)$ and $v\in \mathbb{N}$. The
inverse $\varphi $-transform $T_{\psi }$ is defined by%
\begin{equation*}
T_{\psi }\lambda =\sum_{m\in \mathbb{Z}^{n}}\lambda _{0,m}\Psi
_{m}+\sum_{v=1}^{\infty }\sum_{m\in \mathbb{Z}^{n}}\lambda _{v,m}\psi _{v,m},
\end{equation*}%
where $\lambda =\{\lambda _{v,m}\in \mathbb{C}\}_{v\in \mathbb{N}_{0},m\in 
\mathbb{Z}^{n}}$, see \cite{FJ90}.

To prove of the main result of this paper we need the following $\varphi $%
-transform characterization of $\widetilde{B}_{p(\cdot ),q(\cdot )}^{\alpha
(\cdot ),p(\cdot )}$, see \cite{D6}.

\begin{theorem}
\label{phi-tran}Let $\alpha \in C_{\mathrm{loc}}^{\log }$ and $p,q\in 
\mathcal{P}_{0}^{\log }$ with $0<q^{+}<\infty $. \textit{Suppose that }$\Phi 
$, $\Psi \in \mathcal{S}(\mathbb{R}^{n})$ satisfying $\mathrm{\eqref{Ass1}}$
and $\varphi ,\psi \in \mathcal{S}(\mathbb{R}^{n})$ satisfy $\mathrm{%
\eqref{Ass2}}$ such that $\mathrm{\eqref{Ass4}}$ holds. The operators $%
S_{\varphi }:\widetilde{B}_{p(\cdot ),q(\cdot )}^{\alpha (\cdot ),p(\cdot
)}\rightarrow \widetilde{b}_{p\left( \cdot \right) ,q\left( \cdot \right)
}^{\alpha \left( \cdot \right) ,p(\cdot )}$ and $T_{\psi }:\widetilde{b}%
_{p\left( \cdot \right) ,q\left( \cdot \right) }^{\alpha \left( \cdot
\right) ,p(\cdot )}\rightarrow \widetilde{B}_{p(\cdot ),q(\cdot )}^{\alpha
(\cdot ),p(\cdot )}$ are bounded. Furthermore, $T_{\psi }\circ S_{\varphi }$
is the identity on $\widetilde{B}_{p(\cdot ),q(\cdot )}^{\alpha (\cdot
),p(\cdot )}$.
\end{theorem}

\section{\textbf{D}uality}

This section is devoted to the duality of Triebel-Lizorkin spaces $%
F_{1,q\left( \cdot \right) }^{\alpha \left( \cdot \right) }$. In the case of
constant indices $q$ and $\alpha $, this is a classical part of the theory
of function spaces. Before proving the duality of these function spaces we
present some results, which appeared in the paper of Frazier and Jawerth 
\cite{FJ90} for constant exponents.

\begin{prop}
\label{prop1}\textit{Let }$\alpha \in C_{\mathrm{loc}}^{\log }$\textit{\ , }$%
v\in \mathbb{N}_{0}$, $m\in \mathbb{Z}^{n}$ \textit{and }$p,q\in \mathcal{P}%
_{0}^{\log }$ with $0<q^{-}\leq q^{+}<\infty $\textit{. Suppose that for
each dyadic cube }$Q_{v,m}$ there is a set $E_{Q_{v,m}}\subseteq Q_{v,m}$
with $|E_{Q_{v,m}}|>\varepsilon |Q_{v,m}|$, $\varepsilon >0$. Then%
\begin{equation*}
\left\Vert \lambda \right\Vert _{\widetilde{b}_{p\left( \cdot \right)
,q\left( \cdot \right) }^{\alpha \left( \cdot \right) ,p(\cdot )}}\approx
\sup_{\{P\in \mathcal{Q},|P|\leq 1\}}\Big\|\Big(\frac{\sum\limits_{m\in 
\mathbb{Z}^{n}}2^{v(\alpha \left( \cdot \right) +n/2)}\lambda _{v,m}\chi
_{E_{Q_{v,m}}}}{|P|^{1/p(\cdot )}}\chi _{P}\Big)_{v\geq v_{P}}\Big\|_{\ell
^{q(\cdot )}(L^{p(\cdot )})}.
\end{equation*}
\end{prop}

\textbf{Proof.} Obviously, the problem can be reduced to the case when $\ell
^{q(\cdot )}(L^{p(\cdot )})$ is a normed space. Since $\chi _{E_{Q}}\leq
\chi _{Q}$ for all $Q\in \mathcal{Q}$, one the direction is trivial. For the
other, we use the estimate $\chi _{Q_{v,m}}\leq c$ $\eta _{v,N}\ast \chi
_{E_{Q_{v,m}}}$ for all $Q_{v,m}\in \mathcal{Q}$ and all $N>2n+c_{\log
}(1/p)+c_{\log }(1/q)$. Now Lemma \ref{Alm-Hastolemma1} implies that $%
\left\Vert \lambda \right\Vert _{\widetilde{b}_{p\left( \cdot \right)
,q\left( \cdot \right) }^{\alpha \left( \cdot \right) ,p(\cdot )}}$ is
bounded by%
\begin{eqnarray*}
&&c\sup_{\{P\in \mathcal{Q},|P|\leq 1\}}\Big\|\Big(\frac{\sum\limits_{m\in 
\mathbb{Z}^{n}}2^{v(\alpha \left( \cdot \right) +n/2)}\lambda _{v,m}\eta
_{v,N}\ast \chi _{E_{Q_{v,m}}}}{|P|^{1/p(\cdot )}}\chi _{P}\Big)_{v\geq
v_{P}}\Big\|_{\ell ^{q(\cdot )}(L^{p(\cdot )})} \\
&\lesssim &\sup_{\{P\in \mathcal{Q},|P|\leq 1\}}\Big\|\Big(\frac{%
\sum\limits_{m\in \mathbb{Z}^{n}}2^{v(\alpha \left( \cdot \right)
+n/2)}\lambda _{v,m}\chi _{E_{Q_{v,m}}}}{|P|^{1/p(\cdot )}}\chi _{P}\Big)%
_{v\geq v_{P}}\Big\|_{\ell ^{q(\cdot )}(L^{p(\cdot )})}.
\end{eqnarray*}%
\hspace*{\fill}\rule{3mm}{3mm}

\begin{prop}
\label{Prop2}\textit{Let }$\alpha \in C_{\mathrm{loc}}^{\log }$\textit{\ , }$%
v\in \mathbb{N}_{0}$, $m\in \mathbb{Z}^{n}$ \textit{and }$q\in \mathcal{P}%
_{0}^{\log }$ with $0<q^{-}\leq q^{+}<\infty $\textit{. Suppose that for
each dyadic cube }$Q_{v,m}$ there is a set $E_{Q_{v,m}}\subseteq Q_{v,m}$
with $|E_{Q_{v,m}}|>\varepsilon |Q_{v,m}|$, $\varepsilon >0$. Then%
\begin{equation*}
\left\Vert \lambda \right\Vert _{\widetilde{b}_{q\left( \cdot \right)
,q\left( \cdot \right) }^{\alpha \left( \cdot \right) ,q(\cdot )}}\leq c%
\Big\|\Big(\sum_{v=0}^{\infty }\sum\limits_{m\in \mathbb{Z}^{n}}2^{v(\alpha
\left( \cdot \right) +n/2)q(\cdot )}|\lambda _{v,m}|^{q(\cdot )}\chi
_{E_{Q_{v,m}}}\Big)^{1/q(\cdot )}\Big\|_{\infty }.
\end{equation*}
\end{prop}

\textbf{Proof.} From the fact that $\ell ^{q(\cdot )}(L^{q(\cdot
)})=L^{q(\cdot )}(\ell ^{q(\cdot )})$, we have for any dyadic cube $P$, with 
$|P|\leq 1$%
\begin{eqnarray*}
&&\Big\|\Big(\frac{\sum\limits_{m\in \mathbb{Z}^{n}}2^{v(\alpha \left( \cdot
\right) +n/2)}\lambda _{v,m}\chi _{E_{Q_{v,m}}}}{|P|^{1/q(\cdot )}}\chi _{P}%
\Big)_{v\geq v_{P}}\Big\|_{L^{q(\cdot )}(\ell ^{q(\cdot )})} \\
&=&\Big\|\Big(\frac{1}{|P|}\sum\limits_{v=v_{P}}^{\infty }\sum\limits_{m\in 
\mathbb{Z}^{n}}2^{v(\alpha \left( \cdot \right) +n/2)q(\cdot )}|\lambda
_{v,m}|^{q(\cdot )}\chi _{E_{Q_{v,m}}}\Big)^{1/q(\cdot )}\chi _{P}\Big\|%
_{q(\cdot )}.
\end{eqnarray*}%
Obviously this term is bounded by%
\begin{eqnarray*}
&&\Big\|\Big(\sum\limits_{v=v_{P}}^{\infty }\sum\limits_{m\in \mathbb{Z}%
^{n}}2^{v(\alpha \left( \cdot \right) +n/2)q(\cdot )}|\lambda
_{v,m}|^{q(\cdot )}\chi _{E_{Q_{v,m}}}\Big)^{1/q(\cdot )}\Big\|_{\infty }%
\Big\|\frac{1}{|P|^{1/q(\cdot )}}\chi _{P}\Big\|_{q(\cdot )} \\
&\lesssim &\Big\|\Big(\sum\limits_{v=v_{P}}^{\infty }\sum\limits_{m\in 
\mathbb{Z}^{n}}2^{v(\alpha \left( \cdot \right) +n/2)q(\cdot )}|\lambda
_{v,m}|^{q(\cdot )}\chi _{E_{Q_{v,m}}}\Big)^{1/q(\cdot )}\Big\|_{\infty }.
\end{eqnarray*}%
Therefore, we obtain the desired inequality.\hspace*{\fill}\rule{3mm}{3mm}

For any dyadic cube $P$, with $|P|\leq 1$, we set 
\begin{equation*}
G_{P}^{\alpha \left( \cdot \right) ,q(\cdot )}(\lambda )(x):=\Big(%
\sum\limits_{v=v_{P}}^{\infty }\sum\limits_{m\in \mathbb{Z}^{n}}2^{v(\alpha
\left( x\right) +n/2)q(x)}|\lambda _{v,m}|^{q(x)}\chi _{v,m}(x)\Big)%
^{1/q(x)}.
\end{equation*}%
We put%
\begin{equation}
m_{P}^{\alpha \left( \cdot \right) ,q(\cdot )}(\lambda ):=\inf \Big\{%
\varepsilon :|\{x\in P:G_{P}^{\alpha \left( \cdot \right) ,q(\cdot
)}(\lambda )(x)>\varepsilon \}|<\frac{|P|}{4}\Big\}.  \label{inf-set}
\end{equation}%
We also set%
\begin{equation*}
m^{\alpha \left( \cdot \right) ,q(\cdot )}(\lambda
)(x)=\sup_{P}m_{P}^{\alpha \left( \cdot \right) ,q(\cdot )}(\lambda )\chi
_{P}(x).
\end{equation*}%
Then we obtain.

\begin{prop}
\label{prop2}\textit{Let }$\alpha \in C_{\mathrm{loc}}^{\log },q\in \mathcal{%
P}^{\log }$ with $1<q^{-}\leq q^{+}<\infty $\textit{\ and }$\lambda
=\{\lambda _{v,m}\}_{v\in \mathbb{N}_{0},m\in \mathbb{Z}^{n}}\in \widetilde{b%
}_{q\left( \cdot \right) ,q\left( \cdot \right) }^{\alpha \left( \cdot
\right) ,q(\cdot )}$\textit{. Then}%
\begin{equation*}
\left\Vert \lambda \right\Vert _{\widetilde{b}_{q\left( \cdot \right)
,q\left( \cdot \right) }^{\alpha \left( \cdot \right) ,q(\cdot )}}\approx
\left\Vert m^{\alpha \left( \cdot \right) ,q(\cdot )}(\lambda )\right\Vert
_{\infty }.
\end{equation*}
\end{prop}

\textbf{Proof.} We use the arguments of \cite[Proposition 5.5]{FJ90}. Let $P$
be any dyadic cube, with $|P|\leq 1$. We use the Chebyshev inequality,%
\begin{eqnarray*}
|\{x &\in &P:G_{P}^{\alpha \left( \cdot \right) ,q(\cdot )}(\lambda
)(x)>\varepsilon \}| \\
&\leq &\frac{1}{\varepsilon }\int_{P}G_{P}^{\alpha \left( \cdot \right)
,q(\cdot )}(\lambda )(x)dx \\
&\leq &c\frac{\left\Vert \chi _{P}\right\Vert _{q^{\prime }(\cdot )}}{%
\varepsilon }\left\Vert G_{P}^{\alpha \left( \cdot \right) ,q(\cdot
)}(\lambda )\chi _{P}\right\Vert _{q(\cdot )},
\end{eqnarray*}%
by H\"{o}lder's inequality. Using the properties $\mathrm{\eqref{DHHR}}$, $%
\mathrm{\eqref{DHHR1}}$ and $\mathrm{\eqref{equiv-norm2}}$, to estimate the
last expression by%
\begin{eqnarray*}
c\frac{|P|}{\left\Vert \chi _{P}\right\Vert _{q(\cdot )}\varepsilon }%
\left\Vert G_{P}^{\alpha \left( \cdot \right) ,q(\cdot )}(\lambda )\chi
_{P}\right\Vert _{q(\cdot )} &=&c\frac{|P|}{\varepsilon }\Big\|\frac{%
G_{P}^{\alpha \left( \cdot \right) ,q(\cdot )}(\lambda )}{\left\Vert \chi
_{P}\right\Vert _{q(\cdot )}}\chi _{P}\Big\|_{q(\cdot )} \\
&\leq &c\frac{|P|}{\varepsilon }\Big\|\frac{G_{P}^{\alpha \left( \cdot
\right) ,q(\cdot )}(\lambda )}{|P|^{1/q(\cdot )}}\chi _{P}\Big\|_{q(\cdot )}
\\
&\leq &c\frac{|P|}{\varepsilon }\left\Vert \lambda \right\Vert _{\widetilde{b%
}_{q\left( \cdot \right) ,q\left( \cdot \right) }^{\alpha \left( \cdot
\right) ,q(\cdot )}}.
\end{eqnarray*}%
This term is less than to $\frac{|P|}{4}$ if $\varepsilon >c$ $4\left\Vert
\lambda \right\Vert _{\widetilde{b}_{q\left( \cdot \right) ,q\left( \cdot
\right) }^{\alpha \left( \cdot \right) ,q(\cdot )}}$. Hence, 
\begin{equation*}
\left\Vert m^{\alpha \left( \cdot \right) ,q(\cdot )}(\lambda )\right\Vert
_{\infty }\leq c\left\Vert \lambda \right\Vert _{\widetilde{b}_{q\left(
\cdot \right) ,q\left( \cdot \right) }^{\alpha \left( \cdot \right) ,q(\cdot
)}}.
\end{equation*}%
Now let 
\begin{eqnarray*}
&&j(x) \\
&=&\inf \Big\{j\in \mathbb{N}_{0}:\Big(\sum_{v=j}^{\infty }\sum\limits_{h\in 
\mathbb{Z}^{n}}2^{v(\alpha \left( x\right) +n/2)q(x)}|\lambda
_{v,h}|^{q(x)}\chi _{v,h}(x)\Big)^{1/q(x)}\leq m^{\alpha \left( \cdot
\right) ,q(\cdot )}(\lambda )(x)\Big\}.
\end{eqnarray*}%
and%
\begin{eqnarray*}
E_{v,h} &=&\left\{ x\in Q_{v,h}:2^{-j(x)}\geq l(Q_{v,h})\right\} \\
&=&\Big\{x\in Q_{v,h}:G_{Q_{v,h}}^{\alpha \left( \cdot \right) ,q(\cdot
)}(\lambda )(x)\leq m^{\alpha \left( \cdot \right) ,q(\cdot )}(\lambda )(x)%
\Big\}
\end{eqnarray*}%
for any dyadic cube $Q_{v,h}$, $v\in \mathbb{N}_{0}$ and $h\in \mathbb{Z}%
^{n} $. By $\mathrm{\eqref{inf-set}}$, $|E_{Q_{v,h}}|\geq \frac{3|Q_{v,h}|}{4%
}$, and 
\begin{equation*}
\Big(\sum_{v=0}^{\infty }\sum\limits_{h\in \mathbb{Z}^{n}}2^{v(\alpha \left(
x\right) +n/2)q(x)}|\lambda _{v,h}|^{q(x)}\chi _{E_{v,h}}(x)\Big)%
^{1/q(x)}\leq c\text{ }m^{\alpha \left( \cdot \right) ,q(\cdot )}(\lambda
)(x)
\end{equation*}%
for each $x\in \mathbb{R}^{n}$. Multiplying by $|P|^{-1/q(x)}$, $\chi
_{P}(x) $ ($P$ is a dyadic cube such that $Q_{v,h}\subset P$ and $|P|\leq 1$%
) and then taking the $L^{q(\cdot )}$-norm, we obtain 
\begin{eqnarray*}
&&\Big\|\Big(\sum_{v=v_{P}}^{\infty }\sum\limits_{m\in \mathbb{Z}^{n}}\frac{%
2^{v(\alpha \left( \cdot \right) +n/2)q(\cdot )}|\lambda _{v,h}|^{q(\cdot )}%
}{|P|^{1/q(\cdot )}}\chi _{E_{v,h}}\Big)^{1/q(\cdot )}\chi _{P}\Big\|%
_{q(\cdot )} \\
&\lesssim &\left\Vert m^{\alpha \left( \cdot \right) ,q(\cdot )}(\lambda
)\right\Vert _{\infty }\Big\|\frac{1}{|P|^{1/q(\cdot )}}\chi _{P}\Big\|%
_{q(\cdot )}\lesssim \left\Vert m^{\alpha \left( \cdot \right) ,q(\cdot
)}(\lambda )\right\Vert _{\infty }\text{.}
\end{eqnarray*}%
From the last estimate and Proposition \ref{prop1}, we deduce that 
\begin{equation*}
\left\Vert \lambda \right\Vert _{\widetilde{b}_{q\left( \cdot \right)
,q\left( \cdot \right) }^{\alpha \left( \cdot \right) ,q(\cdot )}}\lesssim
\left\Vert m^{\alpha \left( \cdot \right) ,q(\cdot )}(\lambda )\right\Vert
_{\infty }.
\end{equation*}%
\hspace*{\fill}\hspace*{\fill}.\hspace*{\fill}\rule{3mm}{3mm}

By this proposition and Proposition \ref{prop1}, we obtain another
equivalent norm of $\widetilde{b}_{q\left( \cdot \right) ,q\left( \cdot
\right) }^{\alpha \left( \cdot \right) ,q(\cdot )}$.

\begin{prop}
\label{prop3}\textit{Let }$\alpha \in C_{\mathrm{loc}}^{\log }$\textit{\ and 
}$q\in \mathcal{P}^{\log }$ with $1<q^{-}\leq q^{+}<\infty $\textit{. Then }$%
\lambda =\{\lambda _{v,m}\in \mathbb{C}\}_{v\in \mathbb{N}_{0},m\in \mathbb{Z%
}^{n}}\in \widetilde{b}_{q\left( \cdot \right) ,q\left( \cdot \right)
}^{\alpha \left( \cdot \right) ,q(\cdot )}$\textit{\ if and only if for each
\ dyadic cube }$Q_{v,m}$ there is a subset $E_{Q_{v,m}}\subset Q_{v,m}$ with 
$|E_{Q_{v,m}}|>|Q_{v,m}|/2$ (or any other, fixed, number $0<\varepsilon <1$)
such that%
\begin{equation*}
\Big\|\Big(\sum_{v=0}^{\infty }\sum\limits_{m\in \mathbb{Z}^{n}}2^{v(\alpha
\left( \cdot \right) +n/2)q(\cdot )}|\lambda _{v,m}|^{q(\cdot )}\chi
_{E_{v,m}}(\cdot )\Big)^{1/q(\cdot )}\Big\|_{\infty }<\infty .
\end{equation*}%
Moreover, the infimum of this expression over all such collections $%
\{E_{Q_{v,m}}\}_{v,m}$ is equivalent to $\left\Vert \lambda \right\Vert _{%
\widetilde{b}_{q\left( \cdot \right) ,q\left( \cdot \right) }^{\alpha \left(
\cdot \right) ,q(\cdot )}}$.
\end{prop}

Suppose that $1\leq p\leq \infty $, and $1/p+1/q=1$. In the classical
Lebesgue space,%
\begin{equation}
\left\Vert f\right\Vert _{p}=\sup \left\vert \int f(x)g(x)dx\right\vert ,
\label{Lp-duality}
\end{equation}%
where the supremum is taken over all $g\in L^{q}$ with $\left\Vert
g\right\Vert _{q}\leq 1$. Our aim is to extend this result to $\widetilde{b}%
_{q^{\prime }\left( \cdot \right) ,q^{\prime }\left( \cdot \right)
}^{-\alpha \left( \cdot \right) ,q^{\prime }(\cdot )}$, see \cite{CF13} for
variable Lebesgue spaces. Let $q\in \mathcal{P}$, $\alpha :\mathbb{R}%
^{n}\rightarrow \mathbb{R}$. We define the conjugate norm to $\widetilde{b}%
_{q^{\prime }\left( \cdot \right) ,q^{\prime }\left( \cdot \right)
}^{-\alpha \left( \cdot \right) ,q^{\prime }(\cdot )}$. This is the
functional $\left\Vert \cdot \right\Vert _{\widetilde{b}_{q^{\prime }\left(
\cdot \right) ,q^{\prime }\left( \cdot \right) }^{-\alpha \left( \cdot
\right) ,q^{\prime }(\cdot )}}^{\prime }$ given by 
\begin{equation*}
\left\Vert \lambda \right\Vert _{\widetilde{b}_{q^{\prime }\left( \cdot
\right) ,q^{\prime }\left( \cdot \right) }^{-\alpha \left( \cdot \right)
,q^{\prime }(\cdot )}}^{\prime }=\sup \Big|\int_{P}\frac{1}{|P|}%
\sum\limits_{v=v_{P}}^{\infty }\sum\limits_{m\in \mathbb{Z}^{n}}\lambda
_{v,m}s_{v,m}(x)\chi _{v,m}(x)dx\Big|,\quad \lambda =\{\lambda _{v,m}\in 
\mathbb{C}\}_{v\in \mathbb{N}_{0},m\in \mathbb{Z}^{n}}
\end{equation*}%
where the supremum is taking all dyadic cube $P$, with $|P|\leq 1$ and over
all sequence of functions $s=\{s_{v,m}\}_{v\in \mathbb{N}_{0},m\in \mathbb{Z}%
^{n}}$ such that 
\begin{equation*}
\left\Vert s\right\Vert _{\widetilde{b}_{q\left( \cdot \right) ,q\left(
\cdot \right) }^{\alpha \left( \cdot \right) -n,q(\cdot )}}^{\star
}=\sup_{\{P\in \mathcal{Q},|P|\leq 1\}}\Big\|\Big(\frac{\sum\limits_{m\in 
\mathbb{Z}^{n}}2^{v(\alpha \left( \cdot \right) -n/2)}s_{v,m}\left( \cdot
\right) \chi _{v,m}\left( \cdot \right) }{|P|^{1/q(\cdot )}}\chi _{P}\Big)%
_{v\geq v_{P}}\Big\|_{L^{q(\cdot )}(\ell ^{q(\cdot )})}\leq 1.
\end{equation*}%
Let us start with the following lemma.

\begin{lemma}
\label{ball-daulity-norm}Let $q\in \mathcal{P}$, $1<q^{-}\leq q^{+}<\infty $%
, $\alpha :\mathbb{R}^{n}\rightarrow \mathbb{R}$, $\lambda =\{\lambda
_{v,m}\in \mathbb{C}\}_{v\in \mathbb{N}_{0},m\in \mathbb{Z}^{n}}$ and $%
\mathcal{Q}_{0}=\left\{ Q\in \mathcal{Q}:|Q|\leq 1\right\} $. If $\left\Vert
\lambda \right\Vert _{\widetilde{b}_{q^{\prime }\left( \cdot \right)
,q^{\prime }\left( \cdot \right) }^{-\alpha \left( \cdot \right) ,q^{\prime
}(\cdot )}}^{\prime }\leq 1$ and 
\begin{equation*}
\varrho _{\ell ^{q^{\prime }\left( \cdot \right) }(L^{q^{\prime }\left(
\cdot \right) })}\Big(\Big(\frac{1}{|P|^{1/q^{\prime }\left( \cdot \right) }}%
\sum\limits_{m\in \mathbb{Z}^{n}}2^{v(-\alpha \left( \cdot \right)
+n/2)}|\lambda _{v,m}|\chi _{v,m}\Big)_{v\geq v_{P}}\Big)<\infty
\end{equation*}%
for any $P\in \mathcal{Q}_{0}$, then 
\begin{equation*}
\varrho _{\ell ^{q^{\prime }\left( \cdot \right) }(L^{q^{\prime }\left(
\cdot \right) })}\Big(\Big(\frac{1}{|P|^{1/q^{\prime }\left( \cdot \right) }}%
\sum\limits_{m\in \mathbb{Z}^{n}}2^{v(-\alpha \left( \cdot \right)
+n/2)}|\lambda _{v,m}|\chi _{v,m}\Big)_{v\geq v_{P}}\Big)\leq 1
\end{equation*}%
for any $P\in \mathcal{Q}_{0}.$
\end{lemma}

\textbf{Proof.} Assume, for the sake of contradiction, that 
\begin{equation*}
\varrho _{\ell ^{q^{\prime }\left( \cdot \right) }(L^{q^{\prime }\left(
\cdot \right) })}\Big(\Big(\frac{1}{|P|^{1/q^{\prime }\left( \cdot \right) }}%
\sum\limits_{m\in \mathbb{Z}^{n}}2^{v(-\alpha \left( \cdot \right)
+n/2)}|\lambda _{v,m}|\chi _{v,m}\Big)_{v\geq v_{P}}\Big)>1
\end{equation*}%
for $P\in \mathcal{Q}_{1}\subset \mathcal{Q}_{0}$. Then by the continuity of
the modular there exists $d>1$ such that 
\begin{equation*}
\varrho _{\ell ^{q^{\prime }\left( \cdot \right) }(L^{q^{\prime }\left(
\cdot \right) })}\Big(\Big(\frac{1}{|P|^{1/q^{\prime }\left( \cdot \right) }}%
\sum\limits_{m\in \mathbb{Z}^{n}}2^{v(-\alpha \left( \cdot \right)
+n/2)}\left\vert \frac{\lambda _{v,m}}{d}\right\vert \chi _{v,m}\Big)_{v\geq
v_{P}}\Big)=1.
\end{equation*}%
for $P\in \mathcal{Q}_{1}$. Let $s=\{s_{v,m}\}_{v\in \mathbb{N}_{0},m\in 
\mathbb{Z}^{n}}$ be a sequence of functions defined by%
\begin{equation*}
s_{v,m}(x)=\left\{ 
\begin{array}{ccc}
2^{v(-\alpha \left( x\right) +n/2)q^{\prime }(x)}\left\vert \frac{\lambda
_{v,m}}{d}\right\vert ^{q^{\prime }(x)-1}\chi _{v,m}(x)\text{sgn }\lambda
_{v,m} & \text{if} & Q_{v,m}\subset P\in \mathcal{Q}_{1} \\ 
2^{v(-\alpha \left( x\right) +n/2)q^{\prime }(x)}\left\vert \lambda
_{v,m}\right\vert ^{q^{\prime }(x)-1}\chi _{v,m}(x)\text{sgn }\lambda _{v,m}
& \text{if} & Q_{v,m}\subset P\in \mathcal{Q}_{0}\backslash \mathcal{Q}_{1}.%
\end{array}%
\right. 
\end{equation*}%
Then 
\begin{eqnarray*}
&&\varrho _{\ell ^{q(\cdot )}(L^{q\left( \cdot \right) })}\Big(\Big(\frac{1}{%
|P|^{1/q\left( \cdot \right) }}\sum\limits_{m\in \mathbb{Z}^{n}}2^{v(\alpha
\left( \cdot \right) -n/2)}|s_{v,m}|\Big)_{v\geq v_{P}}\Big) \\
&=&\varrho _{\ell ^{q^{\prime }\left( \cdot \right) }(L^{q^{\prime }\left(
\cdot \right) })}\Big(\Big(\frac{1}{|P|^{1/q^{\prime }\left( \cdot \right) }}%
\sum\limits_{m\in \mathbb{Z}^{n}}2^{v(-\alpha \left( \cdot \right)
+n/2)}\left\vert \frac{\lambda _{v,m}}{d}\right\vert \chi _{v,m}\Big)_{v\geq
v_{P}}\Big) \\
&=&1,
\end{eqnarray*}%
if $P\in \mathcal{Q}_{1}$\ and 
\begin{equation*}
\varrho _{\ell ^{q(\cdot )}(L^{q\left( \cdot \right) })}\Big(\Big(\frac{1}{%
|P|^{1/q\left( \cdot \right) }}\sum\limits_{m\in \mathbb{Z}^{n}}2^{v(\alpha
\left( \cdot \right) -n/2)}|s_{v,m}|\Big)_{v\geq v_{P}}\Big)\leq 1,
\end{equation*}%
if $P\in \mathcal{Q}_{0}\backslash \mathcal{Q}_{1}$, so\ $\left\Vert
s\right\Vert _{\widetilde{b}_{q\left( \cdot \right) ,q\left( \cdot \right)
}^{\alpha \left( \cdot \right) -n,q(\cdot )}}^{\star }\leq 1$. Therefore, by
the definition of $\left\Vert \lambda \right\Vert _{\widetilde{b}_{q^{\prime
}\left( \cdot \right) ,q^{\prime }\left( \cdot \right) }^{-\alpha \left(
\cdot \right) ,q^{\prime }(\cdot )}}^{\prime }$,%
\begin{equation*}
\left\Vert \lambda \right\Vert _{\widetilde{b}_{q^{\prime }\left( \cdot
\right) ,q^{\prime }\left( \cdot \right) }^{-\alpha \left( \cdot \right)
,q^{\prime }(\cdot )}}^{\prime }\geq \Big|\int_{P}\frac{1}{|P|}%
\sum\limits_{v=v_{P}}^{\infty }\sum\limits_{m\in \mathbb{Z}^{n}}\lambda
_{v,m}s_{v,m}(x)\chi _{v,m}(x)dx\Big|,
\end{equation*}%
for any $P\in \mathcal{Q}_{1}$. But the last expression is%
\begin{eqnarray*}
&&\frac{d}{|P|}\sum\limits_{v=v_{P}}^{\infty }\Big\|\sum\limits_{m\in 
\mathbb{Z}^{n}}2^{v(-\alpha \left( \cdot \right) +n/2)q^{\prime }(\cdot
)}\left\vert \frac{\lambda _{v,m}}{d}\right\vert ^{q^{\prime }(\cdot )}\chi
_{v,m}\chi _{P}\Big\|_{1} \\
&=&d\text{ }\varrho _{\ell ^{q^{\prime }\left( \cdot \right) }(L^{q^{\prime
}\left( \cdot \right) })}\Big(\Big(\frac{1}{|P|^{1/q^{\prime }\left( \cdot
\right) }}\sum\limits_{m\in \mathbb{Z}^{n}}2^{v(-\alpha \left( \cdot \right)
+n/2)}\left\vert \frac{\lambda _{v,m}}{d}\right\vert \chi _{v,m}\Big)_{v\geq
v_{P}}\Big) \\
&=&d.
\end{eqnarray*}%
This contradicts our hypothesis on $\lambda =\{\lambda _{v,m}\in \mathbb{C}%
\}_{v\in \mathbb{N}_{0},m\in \mathbb{Z}^{n}}$ , so the desired inequality
holds.\hspace*{\fill}\rule{3mm}{3mm}

The following lemma is the $\widetilde{b}_{q^{\prime }\left( \cdot \right)
,q^{\prime }\left( \cdot \right) }^{-\alpha \left( \cdot \right) ,q^{\prime
}(\cdot )}$-version of $\mathrm{\eqref{Lp-duality}}$, and then we obtain an
equivalent norm on $\widetilde{b}_{q^{\prime }\left( \cdot \right)
,q^{\prime }\left( \cdot \right) }^{-\alpha \left( \cdot \right) ,q^{\prime
}(\cdot )}$.

\begin{lemma}
\label{duality}Let $q\in \mathcal{P}$, $1<q^{-}\leq q^{+}<\infty $ and $%
\alpha :\mathbb{R}^{n}\rightarrow \mathbb{R}$\textit{.} Then%
\begin{equation*}
\left\Vert \lambda \right\Vert _{\widetilde{b}_{q^{\prime }\left( \cdot
\right) ,q^{\prime }\left( \cdot \right) }^{-\alpha \left( \cdot \right)
,q^{\prime }(\cdot )}}\leq \left\Vert \lambda \right\Vert _{\widetilde{b}%
_{q^{\prime }\left( \cdot \right) ,q^{\prime }\left( \cdot \right)
}^{-\alpha \left( \cdot \right) ,q^{\prime }(\cdot )}}^{\prime }\leq
2\left\Vert \lambda \right\Vert _{\widetilde{b}_{q^{\prime }\left( \cdot
\right) ,q^{\prime }\left( \cdot \right) }^{-\alpha \left( \cdot \right)
,q^{\prime }(\cdot )}}.
\end{equation*}
\end{lemma}

\textbf{Proof.}\ Since $\frac{1}{q(\cdot )}+\frac{1}{q^{\prime }\left( \cdot
\right) }=1$, by H\"{o}lder's inequality%
\begin{equation*}
\left\Vert \lambda \right\Vert _{\widetilde{b}_{q^{\prime }\left( \cdot
\right) ,q^{\prime }\left( \cdot \right) }^{-\alpha \left( \cdot \right)
,q^{\prime }(\cdot )}}^{\prime }\leq 2\left\Vert \lambda \right\Vert _{%
\widetilde{b}_{q^{\prime }\left( \cdot \right) ,q^{\prime }\left( \cdot
\right) }^{-\alpha \left( \cdot \right) ,q^{\prime }(\cdot )}}\left\Vert
s\right\Vert _{\widetilde{b}_{q\left( \cdot \right) ,q\left( \cdot \right)
}^{\alpha \left( \cdot \right) -n,q(\cdot )}}^{\star }\leq 2\left\Vert
\lambda \right\Vert _{\widetilde{b}_{q^{\prime }\left( \cdot \right)
,q^{\prime }\left( \cdot \right) }^{-\alpha \left( \cdot \right) ,q^{\prime
}(\cdot )}}\text{.}
\end{equation*}%
Now let us prove that $\left\Vert \lambda \right\Vert _{\widetilde{b}%
_{q^{\prime }\left( \cdot \right) ,q^{\prime }\left( \cdot \right)
}^{-\alpha \left( \cdot \right) ,q^{\prime }(\cdot )}}\leq \left\Vert
\lambda \right\Vert _{\widetilde{b}_{q^{\prime }\left( \cdot \right)
,q^{\prime }\left( \cdot \right) }^{-\alpha \left( \cdot \right) ,q^{\prime
}(\cdot )}}^{\prime }$. By the scaling argument, it suffices to consider the
case $\left\Vert \lambda \right\Vert _{\widetilde{b}_{q^{\prime }\left(
\cdot \right) ,q^{\prime }\left( \cdot \right) }^{-\alpha \left( \cdot
\right) ,q^{\prime }(\cdot )}}^{\prime }\leq 1$ and show that the modular of
the sequence $\lambda $ on the left-hand side is bounded. By Lemma \ref%
{ball-daulity-norm},%
\begin{equation*}
\frac{1}{|P|}\int_{P}\sum\limits_{v=v_{P}}^{\infty }\sum\limits_{m\in 
\mathbb{Z}^{n}}2^{v(-\alpha \left( x\right) +n/2)q^{\prime }(x)}|\lambda
_{v,m}|^{q^{\prime }(x)}\chi _{v,m}(x)dx\leq 1
\end{equation*}%
for any dyadic cube $P$, with $|P|\leq 1$, which is the desired inequality.%
\hspace*{\fill}\rule{3mm}{3mm}

To prove the main result of this paper, we need the following result.

\begin{theorem}
\label{main-result1}\textit{Let }$\alpha \in C_{\mathrm{loc}}^{\log }$ 
\textit{and }$q\in \mathcal{P}^{\log }$ with $1<q^{-}\leq q^{+}<\infty $%
\textit{. Then }%
\begin{equation*}
\big(f_{1,q\left( \cdot \right) }^{\alpha \left( \cdot \right) }\big)^{\ast
}=\widetilde{b}_{q^{\prime }\left( \cdot \right) ,q^{\prime }\left( \cdot
\right) }^{-\alpha \left( \cdot \right) ,q^{\prime }(\cdot )}.
\end{equation*}%
In particular, if $\lambda =\{\lambda _{v,m}\}_{v\in \mathbb{N}_{0},m\in 
\mathbb{Z}^{n}}\in \widetilde{b}_{q^{\prime }\left( \cdot \right) ,q^{\prime
}\left( \cdot \right) }^{-\alpha \left( \cdot \right) ,q^{\prime }(\cdot )}$%
, then the map 
\begin{equation*}
s=\{s_{v,m}\}_{v\in \mathbb{N}_{0},m\in \mathbb{Z}^{n}}\rightarrow
T_{\lambda }(s)=\sum_{v=0}^{\infty }\sum_{m\in \mathbb{Z}^{n}}s_{v,m}\bar{%
\lambda}_{v,m}
\end{equation*}%
defined a continuos linear functional on $f_{1,q\left( \cdot \right)
}^{\alpha \left( \cdot \right) }$ with 
\begin{equation*}
\left\Vert T_{\lambda }\right\Vert _{\big(f_{1,q\left( \cdot \right)
}^{\alpha \left( \cdot \right) }\big)^{\ast }}\approx \left\Vert \lambda
\right\Vert _{\widetilde{b}_{q^{\prime }\left( \cdot \right) ,q^{\prime
}\left( \cdot \right) }^{-\alpha \left( \cdot \right) ,q^{\prime }(\cdot )}},
\end{equation*}%
and every $T\in \big(f_{1,q\left( \cdot \right) }^{\alpha \left( \cdot
\right) }\big)^{\ast }$ is of this form for some $\lambda \in \widetilde{b}%
_{q^{\prime }\left( \cdot \right) ,q^{\prime }\left( \cdot \right)
}^{-\alpha \left( \cdot \right) ,q^{\prime }(\cdot )}$.
\end{theorem}

\textbf{Proof.} We will use some idea from \cite[Theorem 5.9]{FJ90}. Let%
\begin{equation*}
E_{v,h}=\left\{ x\in Q_{v,h}:G_{Q_{v,h}}^{-\alpha \left( \cdot \right)
,q^{\prime }(\cdot )}(\lambda )(x)\leq m^{-\alpha \left( \cdot \right)
,q^{\prime }(\cdot )}(\lambda )(x)\right\}
\end{equation*}%
for any dyadic cube $Q_{v,h}$, with $v\in \mathbb{N}_{0}$ and $h\in \mathbb{Z%
}^{n}$. Then $|E_{Q_{v,h}}|\geq \frac{3|Q_{v,h}|}{4}$ and%
\begin{equation*}
|s_{v,m}||\lambda _{v,h}|=\frac{1}{|E_{v,h}|}\int_{E_{v,h}}|s_{v,h}||\lambda
_{v,h}|dx\leq \frac{4}{3|Q_{v,h}|}\int_{E_{v,h}}|s_{v,h}||\lambda _{v,h}|dx.
\end{equation*}%
Using the H\"{o}lder inequality, we find 
\begin{eqnarray*}
|T_{\lambda }(s)| &\leq &\frac{4}{3}\int \sum_{v=0}^{\infty }\sum_{h\in 
\mathbb{Z}^{n}}2^{v(\alpha \left( x\right) +n/2)}|s_{v,h}|2^{v(-\alpha
\left( x\right) +n/2)}|\lambda _{v,h}|\chi _{E_{v,h}}(x)dx \\
&\leq &\frac{4}{3}\int \Big(\sum_{v=0}^{\infty }\sum_{h\in \mathbb{Z}%
^{n}}2^{v(\alpha \left( x\right) +n/2)q(x)}|s_{v,h}|^{q(x)}\chi _{E_{v,h}}(x)%
\Big)^{1/q(x)} \\
&&\Big(\sum_{v=0}^{\infty }\sum_{h\in \mathbb{Z}^{n}}2^{v(-\alpha \left(
x\right) +n/2)q^{\prime }(x)}|\lambda _{v,h}|^{q^{\prime }(x)}\chi
_{E_{v,h}}(x)\Big)^{1/q^{\prime }(x)}dx.
\end{eqnarray*}%
The last term is bounded by 
\begin{eqnarray*}
&&c\left\Vert s\right\Vert _{f_{1,q\left( \cdot \right) }^{\alpha \left(
\cdot \right) }}\Big\|\Big(\sum_{v=0}^{\infty }\sum_{h\in \mathbb{Z}%
^{n}}2^{v(-\alpha \left( x\right) +n/2)q^{\prime }(\cdot )}|\lambda
_{v,h}|^{q^{\prime }(\cdot )}\chi _{E_{v,h}}\Big)^{1/q^{\prime }(\cdot )}%
\Big\|_{\infty } \\
&\lesssim &\left\Vert s\right\Vert _{f_{1,q\left( \cdot \right) }^{\alpha
\left( \cdot \right) }}\left\Vert m^{-\alpha \left( \cdot \right) ,q^{\prime
}(\cdot )}(\lambda )\right\Vert _{\infty }\lesssim \left\Vert s\right\Vert
_{f_{1,q\left( \cdot \right) }^{\alpha \left( \cdot \right) }}\left\Vert
\lambda \right\Vert _{\widetilde{b}_{q^{\prime }\left( \cdot \right)
,q^{\prime }\left( \cdot \right) }^{-\alpha \left( \cdot \right) ,q^{\prime
}(\cdot )}},
\end{eqnarray*}%
by Proposition \ref{prop2}. Therefore, 
\begin{equation*}
\left\Vert T_{\lambda }\right\Vert _{\big(f_{1,q\left( \cdot \right)
}^{\alpha \left( \cdot \right) }\big)^{\ast }}\lesssim \left\Vert \lambda
\right\Vert _{\widetilde{b}_{q^{\prime }\left( \cdot \right) ,q^{\prime
}\left( \cdot \right) }^{-\alpha \left( \cdot \right) ,q^{\prime }(\cdot )}}.
\end{equation*}%
Clearly every $T\in \big(f_{1,q\left( \cdot \right) }^{\alpha \left( \cdot
\right) }\big)^{\ast }$ is of the form $\sum_{v=0}^{\infty }\sum_{h\in 
\mathbb{Z}^{n}}t_{v,h}\bar{\lambda}_{v,h}$ for some $\lambda =\{\lambda
_{v,h}\}_{v\in \mathbb{N}_{0},h\in \mathbb{Z}^{n}}$. Now, the norm $%
\left\Vert \lambda \right\Vert _{\widetilde{b}_{q^{\prime }\left( \cdot
\right) ,q^{\prime }\left( \cdot \right) }^{-\alpha \left( \cdot \right)
,q^{\prime }(\cdot )}}$ is bounded by%
\begin{equation*}
\sup \Big|\int_{P}\frac{1}{|P|}\sum\limits_{v=v_{P}}^{\infty
}\sum\limits_{h\in \mathbb{Z}^{n}}\lambda _{v,h}s_{v,h}\left( x\right) \chi
_{v,h}(x)dx\Big|,
\end{equation*}%
where the supremum is taking all dyadic cube $P$, with $|P|\leq 1$ and over
all sequence of functions $s=\{s_{v,h}\}_{v\in \mathbb{N}_{0},h\in \mathbb{Z}%
^{n}}$ such that $\left\Vert s\right\Vert _{\widetilde{b}_{q\left( \cdot
\right) ,q\left( \cdot \right) }^{\alpha \left( \cdot \right) -n,q(\cdot
)}}^{\star }\leq 1$, see Lemma \ref{duality}. The integral, $\int_{P}\cdot
\cdot \cdot dx$ can be estimated by%
\begin{equation*}
\sum\limits_{v=v_{P}}^{\infty }\sum\limits_{h\in \mathbb{Z}^{n}}|\lambda
_{v,h}|\int_{P}\frac{1}{|P|}|s_{v,h}(x)|\chi
_{v,h}(x)dx=\sum\limits_{v=0}^{\infty }\sum\limits_{h\in \mathbb{Z}%
^{n}}|\lambda _{v,h}|D_{v,h,P}\leq \left\Vert T\right\Vert _{\big(%
f_{1,q\left( \cdot \right) }^{\alpha \left( \cdot \right) }\big)^{\ast
}}\left\Vert D_{P}\right\Vert _{f_{1,q\left( \cdot \right) }^{\alpha \left(
\cdot \right) }},
\end{equation*}%
where $D_{P}=\{D_{v,h,P}\}_{v\mathbb{\in N}_{0},h\in \mathbb{Z}^{n}}$ and 
\begin{equation*}
D_{v,h,P}=\left\{ 
\begin{array}{ccc}
0 & \text{if} & v<v_{P} \\ 
0 & \text{if} & v\geq v_{P}\text{ and }Q_{v,h}\cap P=\emptyset \\ 
\int_{P}\frac{|s_{v,h}(x)|}{\left\vert P\right\vert }\chi _{v,h}(x)dx & 
\text{if} & v\geq v_{P}\text{ and }Q_{v,h}\subset P.%
\end{array}%
\right.
\end{equation*}%
Let us prove that $\left\Vert D_{P}\right\Vert _{f_{1,q\left( \cdot \right)
}^{\alpha \left( \cdot \right) }}\lesssim 1$. The claim can be reformulated
as showing that%
\begin{equation}
\int_{P}\Big(\sum\limits_{v=v_{P}}^{\infty }\sum\limits_{h\in \mathbb{Z}%
^{n}}2^{v(\alpha \left( y\right) +n/2)q(y)}D_{v,h,P}^{q(y)}\chi _{v,h}(y)%
\Big)^{1/q(y)}dy\lesssim 1.  \label{equa-1}
\end{equation}%
By the fact $2^{v\alpha \left( y\right) }\approx 2^{v\alpha \left( x\right)
} $ for any $x,y\in Q_{v,h}$,%
\begin{eqnarray*}
\frac{2^{v(\alpha \left( y\right) -n/2)}D_{v,h,P}}{\left\vert
Q_{v,h}\right\vert } &\leq &\frac{1}{\left\vert Q_{v,h}\right\vert }%
\int_{Q_{v,h}}\frac{2^{v(\alpha \left( x\right) -n/2)}|s_{v,h}\left(
x\right) |}{\left\vert P\right\vert }\chi _{v,h}(x)dx,\quad y\in
Q_{v,h}\subset P \\
&\lesssim &\eta _{v,N}\ast \Big(\frac{2^{v(\alpha \left( \cdot \right)
-n/2)}|s_{v,h}\left( \cdot \right) |}{\left\vert P\right\vert }\chi _{v,h}%
\Big)(y)
\end{eqnarray*}%
for any $y\in \mathbb{R}^{n}$ and any $N>n$. Therefore, the left-hand side
of $\mathrm{\eqref{equa-1}}$ is bounded by%
\begin{eqnarray*}
&&c\int_{P}\Big(\sum\limits_{v=v_{P}}^{\infty }\Big(\eta _{v,N}\ast
\sum\limits_{h\in \mathbb{Z}^{n}}\frac{2^{v(\alpha \left( \cdot \right)
-n/2)}|s_{v,h}\left( \cdot \right) |}{\left\vert P\right\vert }\chi _{v,h}%
\Big)^{q(y)}\Big)^{1/q(y)}dy \\
&\lesssim &\left\Vert \chi _{P}\right\Vert _{q^{\prime }(\cdot )}\Big\|\Big(%
\sum\limits_{v=v_{P}}^{\infty }\Big(\eta _{v,N}\ast \sum\limits_{h\in 
\mathbb{Z}^{n}}\frac{2^{v(\alpha \left( \cdot \right) -n/2)}|s_{v,h}\left(
\cdot \right) |}{\left\vert P\right\vert }\chi _{v,h}\Big)^{q(\cdot )}\Big)%
^{1/q(\cdot )}\chi _{P}\Big\|_{q(\cdot )} \\
&\lesssim &\left\Vert \chi _{P}\right\Vert _{q^{\prime }(\cdot )}\Big\|\Big(%
\sum\limits_{v=v_{P}}^{\infty }\sum\limits_{h\in \mathbb{Z}^{n}}\Big(%
2^{v(\alpha \left( \cdot \right) -n/2)}\frac{|s_{v,h}\left( \cdot \right) |}{%
|P|}\Big)^{q(\cdot )}\chi _{v,h}\Big)^{1/q(\cdot )}\chi _{P}\Big\|_{q(\cdot
)},
\end{eqnarray*}%
by the H\"{o}lder inequality and Lemma \ref{Alm-Hastolemma2}. We can move $%
\left\Vert \chi _{P}\right\Vert _{q^{\prime }(\cdot )}$ inside the norm and
using the properties $\mathrm{\eqref{DHHR}}$ and $\mathrm{\eqref{DHHR1}}$ to
estimate the last expression by%
\begin{equation*}
c\Big\|\Big(\sum\limits_{v=v_{P}}^{\infty }\sum\limits_{h\in \mathbb{Z}^{n}}%
\Big(2^{v(\alpha \left( \cdot \right) -n/2)}\frac{|s_{v,h}\left( \cdot
\right) |}{\left\Vert \chi _{P}\right\Vert _{q(\cdot )}}\Big)^{q(\cdot
)}\chi _{v,h}\Big)^{1/q(\cdot )}\chi _{P}\Big\|_{q(\cdot )}\lesssim
\left\Vert s\right\Vert _{\widetilde{b}_{q\left( \cdot \right) ,q\left(
\cdot \right) }^{\alpha \left( \cdot \right) -n,q(\cdot )}}^{\star }\leq 1%
\text{.}
\end{equation*}%
Therefore, $\left\Vert \lambda \right\Vert _{\widetilde{b}_{q^{\prime
}\left( \cdot \right) ,q^{\prime }\left( \cdot \right) }^{-\alpha \left(
\cdot \right) ,q^{\prime }(\cdot )}}\lesssim \left\Vert T\right\Vert _{\big(%
f_{1,q\left( \cdot \right) }^{\alpha \left( \cdot \right) }\big)^{\ast }}$
and hence completes the proof of this theorem.\hspace*{\fill}\rule{3mm}{3mm}

Using the notation introduced above, we may now state the main result of
this paper.

\begin{theorem}
\label{main-result2}\textit{Let }$\alpha \in C_{\mathrm{loc}}^{\log }$ 
\textit{and }$q\in \mathcal{P}^{\log }$ with $1<q^{-}\leq q^{+}<\infty $%
\textit{. Then }%
\begin{equation*}
\left( F_{1,q\left( \cdot \right) }^{\alpha \left( \cdot \right) }\right)
^{\ast }=\widetilde{B}_{q^{\prime }\left( \cdot \right) ,q^{\prime }\left(
\cdot \right) }^{-\alpha \left( \cdot \right) ,q^{\prime }(\cdot )}.
\end{equation*}%
In particular, if $g\in \widetilde{B}_{q^{\prime }\left( \cdot \right)
,q^{\prime }\left( \cdot \right) }^{-\alpha \left( \cdot \right) ,q^{\prime
}(\cdot )}$, then the map, given by $l_{g}(f)=\langle f,g\rangle ,$ defined
initially for $f\in \mathcal{S}(\mathbb{R}^{n})$ extends to a continuous
linear functional on $F_{1,q\left( \cdot \right) }^{\alpha \left( \cdot
\right) }$ with $\left\Vert g\right\Vert _{\widetilde{B}_{q^{\prime }\left(
\cdot \right) ,q^{\prime }\left( \cdot \right) }^{-\alpha \left( \cdot
\right) ,q^{\prime }(\cdot )}}\approx \left\Vert l_{g}\right\Vert _{\big(%
F_{1,q\left( \cdot \right) }^{\alpha \left( \cdot \right) }\big)^{\ast }}$
and every $l\in \big(F_{1,q\left( \cdot \right) }^{\alpha \left( \cdot
\right) }\big)^{\ast }$ satisfies $l=l_{g}$ for some $g\in \widetilde{B}%
_{q^{\prime }\left( \cdot \right) ,q^{\prime }\left( \cdot \right)
}^{-\alpha \left( \cdot \right) ,q^{\prime }(\cdot )}$.
\end{theorem}

\textbf{Proof.} Here we use the same arguments of \cite[Theorem 5.13]{FJ90}.
In $\mathrm{\eqref{Ass4}}$ we may choose $\Phi =\Psi \ $and $\varphi =\psi $%
. If $g\in \widetilde{B}_{q^{\prime }\left( \cdot \right) ,q^{\prime }\left(
\cdot \right) }^{-\alpha \left( \cdot \right) ,q^{\prime }(\cdot )}$ and $%
f\in \mathcal{S}(\mathbb{R}^{n})$, then $\langle f,g\rangle =\langle
S_{\varphi }f,S_{\varphi }g\rangle $ and applying Theorems \ref{main-result1}
and \ref{phi-tran} to obtain%
\begin{eqnarray*}
\left\vert \langle f,g\rangle \right\vert &=&\left\vert \langle S_{\varphi
}f,S_{\varphi }g\rangle \right\vert \leq c\left\Vert S_{\varphi
}f\right\Vert _{f_{1,q\left( \cdot \right) }^{\alpha \left( \cdot \right)
}}\left\Vert S_{\varphi }g\right\Vert _{\widetilde{b}_{q^{\prime }\left(
\cdot \right) ,q^{\prime }\left( \cdot \right) }^{-\alpha \left( \cdot
\right) ,q^{\prime }(\cdot )}} \\
&\lesssim &\left\Vert f\right\Vert _{F_{1,q\left( \cdot \right) }^{\alpha
\left( \cdot \right) }}\left\Vert g\right\Vert _{\widetilde{B}_{q^{\prime
}\left( \cdot \right) ,q^{\prime }\left( \cdot \right) }^{-\alpha \left(
\cdot \right) ,q^{\prime }(\cdot )}}.
\end{eqnarray*}%
Hence $\left\Vert l_{g}\right\Vert _{\big(F_{1,q\left( \cdot \right)
}^{\alpha \left( \cdot \right) }\big)^{\ast }}\lesssim \left\Vert
g\right\Vert _{\widetilde{B}_{q^{\prime }\left( \cdot \right) ,q^{\prime
}\left( \cdot \right) }^{-\alpha \left( \cdot \right) ,q^{\prime }(\cdot )}}$%
. Suppose $l\in \big(F_{1,q\left( \cdot \right) }^{\alpha \left( \cdot
\right) }\big)^{\ast }$, then $l_{1}=l\circ T_{\psi }\in \big(f_{1,q\left(
\cdot \right) }^{\alpha \left( \cdot \right) }\big)^{\ast }$, so by Theorem %
\ref{main-result1} there exists $\lambda =\{\lambda _{v,m}\}_{v\in \mathbb{N}%
_{0},m\in \mathbb{Z}^{n}}\in \widetilde{b}_{q^{\prime }\left( \cdot \right)
,q^{\prime }\left( \cdot \right) }^{-\alpha \left( \cdot \right) ,q^{\prime
}(\cdot )}$ such that 
\begin{equation*}
l_{1}(\lambda )=\sum_{v=0}^{\infty }\sum_{m\in \mathbb{Z}^{n}}s_{v,m}\bar{%
\lambda}_{v,m}
\end{equation*}%
for $s=\{s_{v,m}\}_{v\in \mathbb{N}_{0},m\in \mathbb{Z}^{n}}\in f_{1,q\left(
\cdot \right) }^{\alpha \left( \cdot \right) }$ and $\left\Vert \lambda
\right\Vert _{\widetilde{b}_{q^{\prime }\left( \cdot \right) ,q^{\prime
}\left( \cdot \right) }^{-\alpha \left( \cdot \right) ,q^{\prime }(\cdot
)}}\approx \left\Vert l_{1}\right\Vert _{\big(f_{1,q\left( \cdot \right)
}^{\alpha \left( \cdot \right) }\big)^{\ast }}\lesssim \left\Vert
l\right\Vert _{\big(F_{1,q\left( \cdot \right) }^{\alpha \left( \cdot
\right) }\big)^{\ast }}$, since $T_{\psi }$ is bounded. By Theorem \ref%
{phi-tran},%
\begin{equation*}
l_{1}\circ S_{\varphi }=l\circ T_{\psi }\circ S_{\varphi }=l.
\end{equation*}%
Hence putting 
\begin{equation*}
g=T_{\psi }\lambda =\sum_{m\in \mathbb{Z}^{n}}\lambda _{0,m}\Psi
_{m}+\sum_{v=1}^{\infty }\sum_{m\in \mathbb{Z}^{n}}\lambda _{v,m}\psi _{v,m},
\end{equation*}%
we obtain 
\begin{equation*}
l(f)=l_{1}(S_{\varphi }f)=\langle S_{\varphi }f,\lambda \rangle =\langle
f,g\rangle .
\end{equation*}%
Then $l=l_{g}$ and again by Theorem \ref{main-result1},%
\begin{equation*}
\left\Vert g\right\Vert _{\widetilde{B}_{q^{\prime }\left( \cdot \right)
,q^{\prime }\left( \cdot \right) }^{-\alpha \left( \cdot \right) ,q^{\prime
}(\cdot )}}\lesssim \left\Vert \lambda \right\Vert _{\widetilde{b}%
_{q^{\prime }\left( \cdot \right) ,q^{\prime }\left( \cdot \right)
}^{-\alpha \left( \cdot \right) ,q^{\prime }(\cdot )}}\lesssim \left\Vert
l_{g}\right\Vert _{\big(F_{1,q\left( \cdot \right) }^{\alpha \left( \cdot
\right) }\big)^{\ast }}.
\end{equation*}%
\hspace*{\fill}\rule{3mm}{3mm}

\textbf{Acknowledgment}

A great deal of this work has been carried out during the visit of the
author in Jena, Germany. I wish to thank Professor Winfried Sickel for the
hospitality and for the valuable suggestions.

\end{document}